\documentclass[12pt]{article}

\usepackage{amsfonts,amsmath,amsthm,amssymb,amscd}
\usepackage[dvips]{graphicx}

\newtheorem{theorem}{Theorem}
\newtheorem{lemma}[theorem]{Lemma}
\newtheorem{proposition}[theorem]{Proposition}
\newtheorem{corollary}[theorem]{Corollary}

\numberwithin{equation}{section}
\numberwithin{theorem}{section}

\newcommand{\trace}{\qopname\relax o{tr}}

\newcommand{\diag}{\qopname\relax o{diag}}

\newcommand{\bbP}{\mathbb{P}}
\newcommand{\C}{\mathbb{C}}
\newcommand{\Z}{\mathbb{Z}}

\newcommand{\B}{\mathbf{B}}

\newcommand{\p}{\mathbf{p}}
\newcommand{\q}{\mathbf{q}}
\newcommand{\z}{\mathbf{z}}
\newcommand{\m}{\mathbf{m}}
\newcommand{\e}{\mathbf{e}}

\newcommand{\Res}{\textrm{Res}}
\begin{document}

\begin{center}
{\bf 
{\Large 
Isomonodromic Deformations with \\
an Irregular Singularity 
and Hyperelliptic Curves
}
}
\end{center}

\begin{center}
Department of Engineering Science, 
Niihama National College of Technology,\\
7-1 Yakumo-chou, 
Niihama, Ehime, 792-8580. 
\end{center} 

\begin{center}
By 
Kazuhide Matsuda
\end{center}
\vskip 1mm
\par
\quad {\bf Abstract:} \,
In this paper, 
we extend the result of Kitaev and Korotkin \cite{KK} 
to the case where a monodromy-preserving 
deformation has an irregular singularity. 
For the monodromy-preserving deformation, 
we obtain the $\tau$-function 
whose deformation parameters 
are the positions of regular singularities 
and the parameter $t$ of an irregular singularity. 
Furthermore, 
the $\tau$-function is expressed 
by the hyperelliptic $\Theta$ function 
moving the argument $\z$ and the period $\B,$ 
where $t$ and the positions of regular singularities move 
$z$ and $\B,$ respectively.

\section*{Introduction}
In this paper, 
we extend the result of Kitaev and Korotkin \cite{KK} 
to the case where a monodromy-preserving 
deformation has an irregular singularity.  
For the monodromy-preserving deformation, 
we obtain the 
$\tau$-function 
represented 
by the hyperelliptic $\Theta$ function 
moving 
both 
the argument $\z$ 
and 
the period $\B$.  
In \cite{Matsuda}, 
we constructed the $\tau$-function 
by the elliptic $\Theta$ function 
moving the argument $z$ and the period $\Omega$. 
\par
Miwa, Jimbo and Ueno \cite{MM1}
extended 
the work of Schlesinger 
\cite{Sch} 
and 
established a general theory 
of monodromy-preserving deformation 
for a first order matrix system 
of 
ordinary linear differential equations, 
\begin{equation}
\frac{dY}{dx}
=
A(x)Y,
\quad
A(x)
=
\sum_{\nu=1}^{n}
\sum_{k=0}^{r_{\nu}}
\frac{A_{\nu,-k}}{(x-a_{\nu})^{k+1}}
-
\sum_{k=1}^{r_{\infty}}
A_{\infty, -k} x^{k-1},
\end{equation}
having regular or irregular singularities 
of arbitrary rank.
\par
The monodromy data to be preserved are 
\newline
\newline
(i) \quad Stokes multipliers $S_j^{(\nu)} \,\,
(j=1, \ldots , 2r_{\nu}),$
\newline
(ii) \quad connection matrices $C^{(\nu)},$
\newline
(iii) \quad ``exponents of formal monodromy" $T_0^{(\nu)}.$
\newline
\newline
Miwa, Jimbo and Ueno 
found 
a deformation equation 
as a necessary and sufficient condition 
for the monodromy data 
to be independent of  
deformation parameters 
and 
defined 
the $\tau$-function 
for the deformation equation.
\par
Let us explain the relationship 
between 
the $\tau$-function 
and 
the $\Theta$ function. 
Miwa and Jimbo \cite{MM2} 
constructed 
a monodromy-preserving deformation 
with 
irregular singularities 
and 
expressed 
the $\tau$-function 
with 
the $\Theta$ function 
by 
moving 
its argument $\z,$ 
which $t,$ the parameters of the irregular singularities, 
move.
\par 
Kitaev and Korotkin \cite{KK} 
constructed 
a monodromy-preserving deformation of $2\times2$ Fuchsian systems, 
whose deformation parameters are 
the positions of $2g+2$ regular singularities. 
Its $\tau$-function 
is expressed by 
the hyperelliptic $\Theta$ function  
moving the period $\B,$ 
which the positions of regular singularities move. 
When $g=1,$ 
the $\tau$-function 
is 
equivalent to
Picard's solution of 
the sixth Painlev\'e equation 
by the B\"acklund transformation. 
The aim of this paper 
is
to unify Miwa and Jimbo's result 
and Kitaev and Korotkin's result 
in the hyperelliptic case. 
We note that 
Deift, Its, Kapaev and Zhou \cite{DIKZ} 
also 
constructed the isomonodromic deformations 
of $2\times 2$ Fuchsian systems in terms of the hyperelliptic $\Theta$ functions. 
\newline
\par
This paper is organized as follows. 
In Section 1, 
we explain the basic properties of 
the $\Theta$ function, 
the prime-form and the canonical bi-meromorphic differential. 
\par
In Section 2, 
we study 
an ordinary differential equation which is given by
\begin{equation}
\label{equ:diff1}
\frac{d\Psi}{d\lambda}
= 
\left( 
\sum_{j=1}^{2g+2}
\frac{A_{j}}{\lambda-\lambda_{j}}
-B_{-1}
\right)
\Psi(\lambda),
\end{equation}
whose deformation parameters are $\lambda_1,\lambda_2, \ldots \lambda_{2g+2}$ and diagonal elements of $B_{-1}$. 
Then, following Miwa, Jimbo and Ueno \cite{MM1}, 
we introduce the $\tau$-function.
\par
In Section 3,  
we 
solve a class of Riemann-Hilbert 
(inverse monodromy) problem 
for special parameters. 
Furthermore, 
we prove that 
the solution $\Psi(\lambda)$ 
has 
the following special monodromy data: 
\newline
\newline
(i) \quad Stokes multipliers around $\lambda=\infty$ \,\,
$S^{\infty}_1 = S^{\infty}_2 =1,$
\newline
(ii) \quad connection matrices around $\lambda=\lambda_{j}$ \,\,
$C_{j} \,(1 \leq j \leq 2g+2),$
\newline
(iii) \quad the exponents of formal monodromy
$
T_{\infty,0}=0, T_{j,0}=\diag(-\frac14, \frac14) 
\, (1 \leq j \leq 2g+2).
$
\newline
\par 
In Section 4, 
we calculate 
the monodromy-preserving deformation 
which the solution $\Psi(\lambda)$ satisfies 
and compute the coefficients, 
$B_{-1}, A_{j} \,\,(1 \leq j \leq 2g+2).$ 
\par
In Section 5, we find the $\tau$-function. 
Section 5 consists of three subsections. 
Subsection 5.1 is devoted to $H_t$, the Hamiltonian 
at an irregular singular point $\lambda=\infty$. 
Subsection 5.2 
is devoted to Fay's identities and Rauch's variational formulas. 
Subsection 5.3 is devoted to 
$H_{j} \,(1 \leq j \leq 2g+2)$, 
the Hamiltonians on the deformation parameters 
$\lambda_{j} \,\,(1 \leq j \leq 2g+2).$ 
From $H_t$ and $H_j\,\,(1\leq j \leq 2g+2),$ 
we compute the $\tau$-function and prove Theorem 0.1. 
\newline
\par
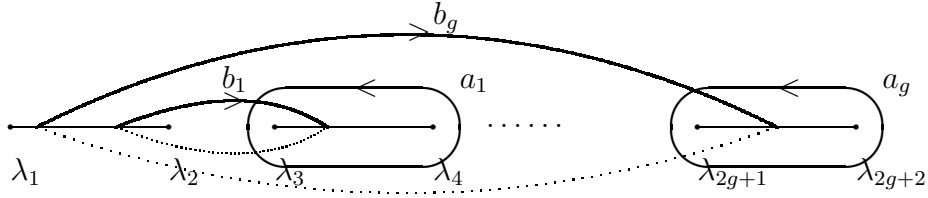
\begin{figure}
\begin{center}
\begin{picture}(400,90)
\put(40,50){\circle*{2}}
\put(40,30){$\lambda_1$ }
\put(100,50){\circle*{2}}
\put(100,30){$\lambda_2$}
\thinlines
\put(40,50){\line(1,0){60}}

\thicklines
\qbezier(80,50)(130,70)(160,50)
\put(120,57){$>$}
\put(120,65){$b_1$}
\thinlines
\qbezier[50](80,50)(130,30)(160,50)

\thicklines
\qbezier(50,50)(190,120)(330,50)
\put(190,83){$>$}

\put(200,89){$b_g$}

\thinlines
\qbezier[80](50,50)(190,0)(330,50)

\put(140,50){\circle*{2}}
\put(140,30){$\lambda_3$ }
\put(200,50){\circle*{2}}
\put(200,30){$\lambda_4$}
\thinlines
\put(140,50){\line(1,0){60}}

\thicklines
\put(170,50){\oval(80,30)}
\put(170,62){$<$}
\put(210,65){$a_1$}

\put(220,50){\ldots \ldots}

\put(300,50){\circle*{2}}
\put(300,30){$\lambda_{2g+1}$}
\put(360,50){\circle*{2}}
\put(360,30){$\lambda_{2g+2}$}
\thinlines
\put(300,50){\line(1,0){60}}

\thicklines
\put(330,50){\oval(80,30)}
\put(330,62){$<$}
\put(370,65){$a_g$}

\end{picture}
\end{center}
\caption
{
Branch cuts 
and 
canonical basis 
of cycles on the hyperelliptic curves, 
$\mathcal{L}$. 
Continuous 
and 
dashed paths 
lie 
on the first and second 
sheet of $\mathcal{L}$, respectively.
}
\label{rs}
\end{figure}

In order to state our main theorem, 
we define  
hyperelliptic curves 
and  
the $\Theta$ function 
and 
introduce 
the prime-form and the canonical bi-meromorphic differential, 
following Fay \cite{Fay}. 
The hyperelliptic curves $\mathcal{L}$ are 
%%%%
\begin{equation*}
\omega^2 = \prod_{j=1}^{2g+2} (\lambda-\lambda_j),
\end{equation*}
whose branch points 
$\lambda_i \in \C \,\,(1 \leq i \leq 2g+2)$ 
are distinct 
and canonical homological basis 
$\{a_j,b_j\}_{1 \leq j \leq g}$ are chosen according to Figure \ref{rs}. 

\par  
The basic holomorphic one forms are expressed by 
\begin{equation*}
dU^0_k = \frac{\lambda^{k-1} d \lambda}{\omega}, \quad 1 \leq k \leq g.
\end{equation*}
Then, the $g \times g$ matrices of a- and b-periods are given by
\begin{equation*}
\mathcal{A}_{kj}= \oint_{a_j} dU^0_k, \quad 
\mathcal{B}_{kj}= \oint_{b_j} dU^0_k, \quad 
1 \leq k,j \leq g.
\end{equation*}
Thus, the normalized holomorphic one forms 
are defined by 
\begin{equation*}
dU_k=\frac{1}{\omega} \sum_{j=1}^g (\mathcal{A}^{-1})_{kj} \lambda^{j-1}d\lambda 
\quad 1 \leq k \leq g,
\end{equation*}
which satisfy 
\begin{equation*}
\oint_{a_j} dU_k = \delta_{jk}.
\end{equation*}
Furthermore, 
from $\mathcal{A}$ and $\mathcal{B}$, 
we can construct 
\begin{equation*}
\bf{B}=\mathcal{A}^{-1} \mathcal{B}.  
\end{equation*}
\par
Following 
Mumford \cite{Mum}, 
we define 
the $\Theta$ function with characteristic $[\p,\q]\,(\p,\q \in \C^g)$ 
by 
\begin{equation}
\label{eqn:theta}
\Theta [\p,\q] (\z  | \B)
=
\sum_{\m \in \Z^g} 
\exp 
\{ 
\pi i 
\langle 
\B (\m+\p),\m+\p 
\rangle 
+ 2 \pi i 
\langle 
\z+\q,\m+\p 
\rangle 
\}, 
\end{equation}
where 
$\z \in \C^g$ and the sum extends over all integer vectors in $\C^g.$
\par
Following Fay \cite{Fay}, 
we introduce the prime-form $E(P,Q)$ 
and the canonical bi-meromorphic differential $W(P,Q)$. 
The prime-form is defined by 
\begin{align*}
E(P,Q)
&=
\frac{\Theta [\mathbf{p}^{*},\mathbf{q}^{*}] (U(P)-U(Q))}{h_{*}(P)h_{*}(Q)}  \\
(h_{*}(P))^2
&=
\sum_{k=1}^g \frac{\partial \Theta [\mathbf{p}^{*},\mathbf{q}^{*}]}{\partial z_k}(0 |\mathbf{B}) {\it dU_k(P)}
\quad
\mathrm{for} \,\, 
P,Q \in \mathcal{L},
\end{align*}
where 
$[\mathbf{p}^{*},\mathbf{q}^{*}]$ is an arbitrary odd non-singular half integer characteristic. 
The canonical bi-meromorphic differential $W(P,Q)$ is given by
\begin{equation*}
W(P,Q)
=
dx_Pdx_Q
\log E(P,Q),
\end{equation*}
where $dx_P,dx_Q$ are the exterior differentiations with respect to the local parameter of $P,Q,$ respectively. 
When $P=\infty^1$ and $Q=\infty^2,$ 
we define the local coordinates of $P,Q$ by 
$$
x_{\infty^1}=\frac{1}{\lambda}, \,\,
x_{\infty^2}=\frac{1}{\lambda}, 
$$ 
and have 
$$
W(\infty^1,\infty^2)
=
\left(
\left.
\frac{\partial^2}{\partial x_{\infty^1}\partial x_{\infty^2}} 
\log E(P,Q) 
\right|_{P=\infty^1,Q=\infty^2}
\right)dx_{\infty^1}dx_{\infty^2}.
$$
\par
According to Fay \cite{Fay}, 
the projective connection $S(Q)$ 
is given by the equation
\begin{equation*}
W(P,Q)
=
\left(
\frac{1}{(x_{P}-x_{Q})^2}
+
\frac{1}{6}
S(Q)
+
O(x_{P}-x_{Q})
\right)
dx_{P}dx_{Q},
\end{equation*}
where $P$ and $Q$ have local coordinates $x_{P},x_{Q}$ 
in a neighborhood of $Q\in \mathcal{L}$. 
When $Q=\infty^1$, 
we define the local coordinate of $Q$ by 
$$
x_{P}=x_{Q}=\frac{1}{\lambda} 
\quad
\mathrm{if} \,\, P \longrightarrow Q=\infty^1.
$$ 
\par
Our main theorem is as follows: 
\begin{theorem}
For the monodromy-preserving deformation (\ref{eqn:mpd}),
the $\tau$-function is
\begin{alignat*}{3}
\tau(\lambda_1, \ldots, \lambda_{2g+2} ; t)
&=
\Theta[\p,\q]
\left(
\mathbf{v}(t) |\B
\right) & 
&(\det \mathcal A)^{-\frac12}
\prod_{1 \leq j<k \leq 2g+2} 
(\lambda_j-\lambda_{k})^{-\frac18} \\
&       &
&\times 
\exp \Big\{\frac{t^2}{4} \left(\frac16 S(\infty^1)-\frac{W(\infty^1,\infty^2)}{dx_{\infty^1}dx_{\infty^2}}\right) \Big\},
\end{alignat*}
where
\begin{equation*}
\mathbf{v}(t)
=
t \times 
\left(
\frac{dU_1}{dx_{\infty^1}}(\infty^1),
\frac{dU_2}{dx_{\infty^1}}(\infty^1),
\ldots,
\frac{dU_{g}}{dx_{\infty^1}}(\infty^1)
\right)
\end{equation*}
and $x_{\infty^1}=\displaystyle\frac{1}{\lambda},$ which is a local coordinate of $\infty^1.$
\end{theorem}

\quad {\bf Remark.}
\,
By setting $t=0$, we obtain Kitaev and Korotkin's $\tau$-function 
in \cite{KK}. We had a branch point at $\infty$ in \cite{Matsuda}, 
while we do not in this paper.

\section{Hyperelliptic Curves and the $\Theta$ Function}
In this section, 
we explain the detailed properties of 
the $\Theta$ function, 
the prime-form $E(P,Q)$ and the canonical bi-meromorphic differential $W(P,Q)$. 
\par
The $\Theta$ function possesses the following periodicity properties:
\begin{equation}
\label{eqn:per1}
\Theta [\p,\q] (\z+\e_j | \B)
=
\exp \{ 2 \pi i p_j \} 
\Theta [\p,\q] (\z | \B)  
\end{equation}

\begin{equation}
\label{eqn:per2}
\Theta [\p,\q] (\z+ \B \e_j | \B)
=
\exp \{ -2 \pi i q_j - \pi i \B_{jj} -2 \pi i z_j \} 
\Theta [\p,\q] (\z | \B),  
\end{equation}
where 
\begin{equation*}
\mathbf{e}_j ={}^t(0, \ldots, \stackrel{j th}{1}, \ldots, 0).
\end{equation*}
\par
We define the Abel map $U(P) \in \C^g \,P\in \mathcal{L}$ by 
\begin{align*}
U(P)
&=
{}^t
\left(
U_1(P), U_2(P), \ldots, U_g(P)
\right) \\
U_j(P)
&=
\int_{\lambda_1}^P dU_j \quad (1\leq j \leq g).
\end{align*}
For the canonical homological basis $\{a_j, b_j \}_{1 \leq j \leq g}$ and the base point $\lambda_1$, 
the Riemann constants $\mathbf{K} \in \C^g$ are as follows:
\begin{equation*}
\mathbf{K}=
\frac12 
\B
\left(
\mathbf{e}_1+\mathbf{e}_2+\cdots+\mathbf{e}_g
\right)
+
\frac12 
\left(
\mathbf{e}_1+2\mathbf{e}_2+\cdots+g\mathbf{e}_g
\right).
\end{equation*}
\par
A characteristic $[\mathbf{p},\mathbf{q}]$ 
is a $g \times 2$ matrix of complex numbers 
which is given by
\begin{equation*}
[\mathbf{p},\mathbf{q}]
=
\begin{bmatrix}
p_1 & q_1 \\
p_2 & q_2 \\
\vdots & \vdots \\
p_g & q_g
\end{bmatrix},
\end{equation*}
where 
\begin{equation*}
\mathbf{p}
=
{}^t
(p_1,p_2, \ldots, p_g),  
\mathbf{q}
=
{}^t
(q_1,q_2,\ldots,q_g).
\end{equation*}
We consider $p_i,q_i \,\,(1\leq i \leq g)$ as elements of 
$\C^g/\Z^g.$ 
If all the components are half-integers, 
$[\p,\q]$ is called a half-integer characteristic. 
A half-integer characteristic is in one-to-one correspondence 
with a half-period $\mathbf{B} \mathbf{p} + \mathbf{q}$. 
If the scalar product $4 \langle \mathbf{p}, \mathbf{q} \rangle$ is odd, 
then the characteristic is called odd 
and 
the related $\Theta$ function is odd with respect to its argument $\mathbf{z}$. 
If this scalar product is even, 
then the characteristic is called even 
and 
the related $\Theta$ function is even with respect to its argument $\mathbf{z}$. 
\par
The odd characteristics which are important for us in the sequel correspond 
to any subset 
\begin{equation*}
S=
\{ \lambda_{i_1}, \lambda_{i_2},\ldots, \lambda_{i_{g-1}}
\}, 
\end{equation*} 
whose components are arbitrary distinct branch points. 
The odd half-period corresponding to the subset $S$ 
is expressed by
\begin{equation*}
\mathbf{B} \mathbf{p}^S + \mathbf{q}^S 
=
U(\lambda_{{\it i_1}})+U(\lambda_{{\it i_2}})+ \cdots +U(\lambda_{{\it i_{g-1}}})-\mathbf{K}.
\end{equation*}
\par
Analogously, 
we shall be interested in the even half-period 
corresponding to the subset 
$$
T=
\{
\lambda_{i_1}, \lambda_{i_2},\ldots, \lambda_{i_{g+1}}
\},
$$  
which consists of arbitrary $g + 1$ branch points. 
The even half-period 
is given by 
\begin{equation*}
\mathbf{B} \mathbf{p}^T + \mathbf{q}^T 
=
U(\lambda_{{\it i_1}})+U(\lambda_{{\it i_2}})+ \cdots +U(\lambda_{{\it i_{g+1}}})-\mathbf{K}.
\end{equation*}
We fix the choice of $T$ and define 
\begin{equation*}
\{\lambda_{j_1},\lambda_{j_2},\ldots,\lambda_{j_{g+1}} \}
:=
\{\lambda_1,\lambda_2,\ldots,\lambda_{2g+2}\}
\setminus T.
\end{equation*} 
\par
We explain the detailed property of $E(P,Q)$. 
On page 13--14 of Fay \cite{Fay}, 
we find
\begin{align*}
E(P,Q)
&=
\frac{ \Theta [\bf{p}^T,\bf{q}^T]({\it U(P)-U(Q)})}{\Theta [\bf{p}^T,\bf{q}^T](0) {\it m_T(P,Q)}},  \\
m_{T}(P,Q)
&=
\frac{
\omega(Q) \prod_{k=1}^{g+1}(\lambda(P)-\lambda_{i_k})
+
\omega(P) \prod_{k=1}^{g+1}(\lambda(Q)-\lambda_{i_k})
}
{2 (\lambda(Q)-\lambda(P))} \\
& 
\hspace{20mm}
\times
\left[
\frac{
d \lambda(P) d \lambda(Q)
}
{
\omega(P)\omega(Q) 
\prod_{k=1}^{g+1}(\lambda(P)-\lambda_{i_k})(\lambda(Q)-\lambda_{i_k})
}
\right]^{\frac12},
\end{align*}
because $\mathcal{L}$ is hyperelliptic.

\section{The Schlesinger System}
We 
study an ordinary differential equation 
which is given by
\begin{equation}
\label{equ:diff2}
\frac{d \Psi}{d \lambda}
=
\left( 
\sum_{j=1}^{2g+2}
\frac{A_j}{\lambda-\lambda_{j}}
-B_{-1}
\right)
\Psi(\lambda),
\end{equation}
where $A_1,A_2,\ldots,A_{2g+2}, B_{-1} \in sl(2, \C)$ 
are independent of $\lambda$. 
\par
The monodromy data of \eqref{equ:diff2} 
are as follows:
\newline
\newline
(i) \quad Stokes multipliers around $\lambda=\infty$ \,\, 
$S^{\infty}_1 = S^{\infty}_2 ;$
\newline
(ii) \quad connection matrices around $\lambda_j$ \,\,
$C_j,  \,(j=1,2,\ldots, 2g+2);$
\newline
(iii) \quad the exponents of formal monodromy 
$T_{\infty,0}, T_{j,0} \, (j=1,2,\ldots, 2g+2).$
\newline
\newline
In the next section, 
we obtain a 
convergent series around 
$\lambda = \infty,  \lambda_{j} \,\,(j=1,2,\ldots, 2g+2)$ 
which are expressed by
\begin{align}
& 
\Psi(\lambda) 
=
(
1
+
O(\frac{1}{\lambda})
)
\exp
T^{\infty}(\lambda)
=\hat{\Psi}^{\infty}(\lambda) 
\exp T^{\infty}(\lambda), \\
& 
\Psi(\lambda)
=
G_{j}
(
1
+
O(\lambda-\lambda_{j})
)
\exp T_{j}(\lambda) 
=
G_{j}
\hat{\Psi}_{j}(\lambda)
\exp
T_{j,0}(\lambda),  \notag \\
&  
\hspace{90mm}(j=1,2,\ldots, 2g+2),
\end{align}
where
\begin{align}
T^{\infty}(\lambda)
&=
\left(
\begin{array}{cc}
-\frac{t}{2} &                   \\
                  & \frac{t}{2}  \\
\end{array}
\right)
\lambda
+
T_{\infty,0}
\log(\frac{1}{\lambda}),  \\
T_{j}(\lambda)
&=
T_{j,0}
\log(\lambda-\lambda_{j})
\,\, (j=1,2,\ldots, 2g+2).
\end{align}
\par
For the deformation parameters 
$t, \lambda_1, \lambda_2, \ldots, \lambda_{2g+2},$ 
the closed one form is defined by
\begin{align}
\Omega 
&=
\omega_{\infty} 
+
\omega_{\lambda_1}
+
\omega_{\lambda_2}
+
\cdots
+
\omega_{\lambda_{2g+2}} \\
&= 
H_t dt 
+ 
H_1 d \lambda_1 
+ 
H_2 d \lambda_2 
+
\cdots
+ 
H_{2g+2} d \lambda_{2g+2}, 
\end{align}
where
\begin{align}
\omega_{\infty}
&=
- 
\Res_{\lambda=\infty} \,\, 
\trace
\hat{\Psi}^{\infty}(\lambda)^{-1} 
\frac{\partial \hat{\Psi}^{\lambda}}{\partial \lambda}
(\lambda)
\, 
d T^{\infty}(\lambda),     \\
\omega_{\lambda_j}
&= 
- 
\Res_{\lambda=\lambda_{j}} \,\, 
\trace
\hat{\Psi}_{j}(\lambda)^{-1} 
\frac{\partial \hat{\Psi}_{j}}{\partial \lambda}(\lambda)
\, 
d T_{j}(\lambda) \,\,(j=1,2,\ldots, 2g+2 ),
\end{align}
and
$d$ 
is the exterior differentiation 
with respect to 
the deformation parameters $t, \lambda_1, \lambda_2, \ldots, \lambda_{2g+2}$. 
Especially, we can write
\begin{equation}
\omega_{\lambda_{j}} 
= 
\left[ 
\Res_{\lambda=\lambda_{j}} 
\frac12 
\trace 
\left( 
\frac{d \Psi}{d \lambda}\Psi^{-1} 
\right)^2 
\right]
d \lambda_{j}.
\end{equation}
Then, from the closed one form $\Omega$, 
the $\tau$-function is defined by
\begin{equation}
\Omega := d \log \tau (\lambda_1, \lambda_2, \ldots, \lambda_{2g+2},t).
\end{equation}

\section{The Riemann-Hilbert Problem for Special Parameters}
In this section, 
we concretely construct 
a $2\times2$ matrix valued function $\Psi(\lambda)$, 
whose monodromy data, 
(i)\, Stokes multipliers, 
(ii)\, connection matrices, 
(iii)\, exponents of formal monodromy, 
are independent of 
the deformation parameters $t, \lambda_1, \lambda_2, \ldots, \lambda_{2g+2}$.

\begin{figure}
\begin{center}
\begin{picture}(400,150)

\put(70,70){\line(1,1){120}}
\put(100,70){\line(3,4){90}}

\put(85,70){\oval(30,55)[b]}
\put(85,60){\circle*{2}}
\put(80,50){$\lambda_1$}
\put(81,39){$>$}
\put(82,25){$l_{1}$}

\put(130,70){\line(1,2){60}}
\put(160,70){\line(1,4){30}}

\put(145,70){\oval(30,55)[b]}
\put(145,60){\circle*{2}}
\put(142,50){$\lambda_2$}
\put(139,39){$>$}
\put(142,25){$l_2$}

\put(180,39){$\cdots$}

\put(190,190){\circle*{2}}
\put(190,195){$\lambda_0$}

\put(220,70){\line(-1,4){30}}
\put(250,70){\line(-1,2){60}}

\put(235,70){\oval(30,55)[b]}
\put(235,60){\circle*{2}}
\put(232,50){$\lambda_{2g+1}$}
\put(232,39){$>$}
\put(232,25){$l_{2g+1}$}

\put(280,70){\line(-3,4){90}}
\put(310,70){\line(-1,1){120}}

\put(295,70){\oval(30,55)[b]}
\put(295,60){\circle*{2}}
\put(292,50){$\lambda_{2g+2}$}
\put(292,39){$>$}
\put(292,25){$l_{2g+2}$}

\end{picture}
\end{center}
\caption
{
Generators of 
$
\pi_1 
\left(
\C \bbP^1
\setminus
\{
\lambda_1,
\ldots,
\lambda_{2g+2}
\},
\lambda_0
\right)
$
}
\end{figure}

The involution of $\mathcal{L}$ is defined by
\begin{equation*}
*:
(\lambda,\omega)
\longrightarrow
(\lambda,-\omega).
\end{equation*}
Then, 
from the $\Theta$ function and $W(P,Q)$, 
we define the $2 \times 2$ matrix valued function 
$\Phi(P)$ 
by
\begin{equation}
\Phi (P) = 
\left(
\begin{array}{cc}
\varphi(P) & \varphi(P^{*}) \\
\psi(P) & \psi(P^{*})
\end{array}
\right),
\end{equation}
where
\begin{align*}
\varphi(P) 
&=
\Theta [\p,\q] (U(P)+U(P_{\varphi})+\mathbf{v}(t)) 
\Theta [\p^S,\q^S](U(P)-U(P_{\varphi})) \\
& 
\quad
\times
\exp 
-\frac{t}{2} 
\{ 
-
\int^P_{\lambda_1}
W(P,\infty^1)
+
\int^P_{\lambda_1}
W(P,\infty^2)
\},  \\
\psi(P) 
&= 
\Theta [\p,\q] (U(P)+U(P_{\psi})+\mathbf{v}(t)) 
\Theta [\p^S, \q^S] (U(P)-U(P_{\psi})) \\
& 
\quad
\times
\exp 
-\frac{t}{2} 
\{ 
-
\int^P_{\lambda_1}
W(P,\infty^1)
+
\int^P_{\lambda_1}
W(P,\infty^2)
\}, 
\\
\mathbf{v}(t)
&=
t \times 
\left(
\frac{dU_1}{dx_{\infty^1}} (\infty^1),
\frac{dU_2}{dx_{\infty^1}} (\infty^1),
\ldots,
\frac{dU_{g}}{dx_{\infty^1}} (\infty^1)
\right),
\end{align*}
where $P_{\varphi}, P_{\psi}$ are arbitrary points of $\mathcal{L}$ 
and $x_{\infty^1}=\displaystyle\frac{1}{\lambda},$ 
which is a local coordinate of $\infty^1.$

\begin{proposition}
\label{prop:reg}
(1) \quad 
The function $\Phi(P)$ is 
invertible outside of 
the branch points $\lambda_1, \lambda_2, \ldots, \lambda_{2g+2}.$ 
\newline
(2) \quad  
$\det \Phi(P)$ has zeros at $\lambda_j \notin S$ 
with the first order 
and 
has zeros at $\lambda_j \in S$ with the third order.
\newline
(3) 
\quad 
The function $\Phi(P)$ 
transforms as follows with respect to 
the tracing along the canonical homological basis, 
$a_j, b_j\,(j=1,2,\ldots, g)$:
\begin{align}
\label{eqn:trans1}
T_{a_j} 
\big[ 
\Phi (P) 
\big] 
&=
\Phi (P) 
\left(
\begin{array}{cc}
\exp \{ 2\pi i (p_j+p^S_j) \} &        \\ 
           & \exp \{ - 2\pi i (p_j+p^S_j) \} 
\end{array}
\right),  \\
\label{eqn:trans2}
T_{b_j} 
\big[ 
\Phi (P) 
\big] 
&=
\Phi (P) 
\left(
\begin{array}{cc}
\exp \{ - 2\pi i (q_j+q^S_j) \} &         \\
              &  \exp \{ 2\pi i (q_j+q^S_j) \}
\end{array}
\right)
\exp \{ -2 \pi i \B_{jj}-4 \pi i U_j(P)\},
\end{align} 
where $T_{l}$ denotes 
the operator of analytic continuation 
along the contour $l$.

\end{proposition}

\begin{proof}
By using the periodicity 
(\ref{eqn:per1}) 
and 
(\ref{eqn:per2}), 
we obtain 
\begin{align*}
T_{a_j}[\varphi(P)]&=
\exp\{2\pi i (p_j+p_j^S)\} \varphi(P),  \\
T_{b_j}[\varphi(P)]&=
\exp\{-2\pi i (q_j+q_j^S)-2 \pi i \B_{jj}-4 \pi i U_j(P)\} \varphi(P).
\end{align*}
We deduce the same transformation laws for $\psi(P)$. 
\par
The actions of the involution $*$ on $\{a_j,b_j\}$ and $dU_j$ 
are given by
\begin{equation*}
a_j^{*}=-a_j,\,
b_j^{*}=-b_j,\,
dU_j(P^{*})=-dU_j(P),\, 
(j=1,2,\ldots,g),
\end{equation*}
respectively.
Therefore, 
we get
\begin{align}
T_{a_j}[\varphi(P^*)]
&=
\exp \{-2\pi i(p_j+p_j^S)\}
\varphi(P^*), \\
T_{b_j}[\varphi(P^*)]&= \exp\{2\pi i(q_j+q_j^S)-2\pi i \B_{jj}-4\pi i U_j(P)\}\varphi(P^*).
\end{align}
We deduce the same transformation laws for $\psi(P^*)$. 
Then, we complete the proof of (3).
\par
The equations (\ref{eqn:trans1}) and (\ref{eqn:trans2}) 
imply that
\begin{align*}
T_{a_j}[\det \Phi(P)]&=\det \Phi(P),\\ 
T_{b_j}[\det \Phi(P)]&=\det \Phi(P) 
\exp \{-4\pi i \B_{jj}-8\pi i U_j(P)\},
\end{align*}
which imply that 
\begin{equation*}
\frac{1}{2 \pi i}
\oint_{\partial \mathcal{L}}
\frac
{
d(\det \Phi(P))
}
{\det \Phi(P)}
=
4g.
\end{equation*}
Thus, 
it follows that 
\begin{equation*}
3(g-1)+g+3=4g,
\end{equation*}
because $\det \Phi(P)$ has zeros at the branch points $\lambda_j$ 
and has zeros at $\lambda_j \in S$ of the order of at least three. 
Therefore, 
it follows that 
$\det \Phi(P)$ does not vanish outside of the branch points 
and 
that 
$\det \Phi(P)$ has zeros at $\lambda_j \notin S$ 
with the first order 
and 
has zeros at $\lambda_j \in S$ with the third order, 
which complete the proof of (1) and (2), respectively.
\end{proof}

In order to normalize $\Phi(\lambda)$ near $\lambda=\infty$, 
we have

\begin{lemma}
\label{lem:Fay1}
For $P,Q \in \mathcal{L},$
\begin{align}
W(P,Q)
&=
\left(
\frac{1}{(x_{P}-x_{Q})^2}
+
\frac{1}{6}
S(Q)
+
O(x_{P}-x_{Q})
\right)
dx_{P}dx_{Q}, \\
\label{eqn:proj2}
\frac16 S(Q)
&=
\frac16 \{\lambda,x_{Q}\}(Q)
+
\frac{1}{16}
\left(
\frac{d}{dx_{Q}}
\log 
\prod_{k=1}^{g+1}
\frac{\lambda-\lambda_{i_k}}{\lambda-\lambda_{j_k}}(Q) 
\right)^2   \notag  \\
& \hspace{30mm}
-\sum_{k,l=1}^g 
\frac{\partial^2}{\partial z_k \partial z_l}
\log \Theta [\p^T,\q^T](0) 
\frac{dU_k}{dx_{Q}}(Q)\frac{dU_l}{dx_{Q}}(Q),
\end{align}
where $x_{P},x_{Q}$ are local coordinates of $P,Q \in \mathcal{L}$ 
and 
\begin{equation*}
\{\lambda,x\}
=
\frac{\lambda^{\prime \prime \prime}}{\lambda^{\prime}}
-\frac32
\left(
\frac{\lambda^{\prime \prime}}{\lambda^{\prime}}
\right)^2.
\end{equation*} 
\end{lemma}

\begin{proof}
See pp. 20 in \cite{Fay}.
\end{proof}
In Lemma \ref{lem:Fay1}. we set 
\begin{equation*}
x_{P}=x_{Q}=
\frac{1}{\lambda},
\quad
Q=\infty^1,
\end{equation*}
and take the limit $P \longrightarrow \infty^1.$ 
Furthermore, we define the constant terms $c_{\infty^1},c_{\infty^2}$ by
\begin{alignat*}{11}
& 
\int_{\lambda_1}^P
W(P,\infty^1)
& 
&=
-\lambda & 
+&c_{\infty^1} &
&+\frac16 S(\infty^1)
\lambda^{-1} & 
&+\cdots & \quad &\mathrm{near} \,\,\lambda=\infty^1, \\
& 
\int_{\lambda_1}^P
W(P,\infty^2)
& 
&=
         &    
&c_{\infty^2} & 
&+
\frac{
W(\infty^1,\infty^2)
}{dx_{\infty^1}dx_{\infty^2}}\lambda^{-1}  &  
&+\cdots  & &\mathrm{near} \,\,\lambda=\infty^1. 
\end{alignat*}
Therefore, 
$\Phi(\lambda)$ can be developed near $\lambda=\infty$ as

\begin{equation*}
\Phi(\lambda)
=
\left(
G^{\infty}
+
O(\frac{1}{\lambda})
\right)
\exp
\left(
\begin{array}{cc}
- 
\frac{t}{2} \lambda&     \\
                    & \frac{t}{2}\lambda  \\
\end{array}
\right),
\end{equation*}
where $G^{\infty}$ is 
a $2 \times 2$ matrix whose matrix elements are given by
\begin{align*}
(G^{\infty})_{11}
&=
\Theta[\p,\q](U(\infty^1)+U(P_{\varphi})+\mathbf{v}(t))
\Theta[\p^S,\q^S](U(\infty^1)-U(P_{\varphi}))
\exp 
\Big\{
\frac{t}{2}
\left(
c_{\infty^1}-c_{\infty^2}
\right)
\Big\},  \\
(G^{\infty})_{21}
&= 
\Theta[\p,\q](U(\infty^1) + U(P_{\psi}) +\mathbf{v}(t)) 
\Theta[\p^S,\q^S](U(\infty^1)-U(P_{\psi})) 
\exp 
\Big\{
\frac{t}{2}
\left(
c_{\infty^1}-c_{\infty^2}
\right) 
\Big\}, \\
(G^{\infty})_{12}
&=
\Theta[\p,\q](U(\infty^2) +U(P_{\varphi}) +\mathbf{v}(t)) 
\Theta[\p^S,\q^S](U(\infty^2) -U(P_{\varphi}))
\exp 
\Big\{
-
\frac{t}{2}
\left(
c_{\infty^1}-c_{\infty^2}
\right)
\Big\},   \\
(G^{\infty})_{22}
&=
\Theta[\p,\q](U(\infty^2) +U(P_{\psi}) +\mathbf{v}(t)) 
\Theta[\p^S,\q^S](U(\infty^2)-U(P_{\psi}))
\exp 
\Big\{
-
\frac{t}{2}
\left(
c_{\infty^1}-c_{\infty^2}
\right)
\Big\}.
\end{align*}
Proposition \ref{prop:reg} shows that 
\begin{equation*}
\det G^{\infty}
=
\det \Phi (\infty)
\neq
0.
\end{equation*}
We define a matrix valued function $\Psi(\lambda)$ 
by
\begin{equation}
\Psi(P)
=
\frac
{
\sqrt{\det \Phi(\infty)}
}
{
\sqrt{\det \Phi(P)}
} 
(G^{\infty})^{-1} 
\Phi(P).
\end{equation}

The expansions of $\Psi(\lambda)$ near $\lambda=\infty$ 
are as follows:

\begin{lemma}
\label{lem:ya}
\begin{align*}
& 
\Psi(\lambda)
=
\left(
\left(
\begin{array}{cc}
1 & 0 \\
0 & 1   
\end{array}
\right)
+
O(\frac{1}{\lambda})
\right)
\exp 
T^{\infty}(\lambda), \\
& 
T^{\infty}(\lambda)
=
T^{\infty}_{-1}
\lambda, \,\, 
T^{\infty}_{-1}
=
\left(
\begin{array}{cc}
-\frac{t}{2} &      \\
                 & \frac{t}{2}
\end{array}
\right), 
\end{align*}
where the Taylor series 
of 
$\Psi(\lambda) \exp\left\{-T^{\infty}(\lambda)\right\}$ 
is convergent. 
Especially, if 
$P_{\varphi}=\infty^1$ and $ P_{\psi}=\infty^2,$

\begin{align*}
\Psi(\lambda)
&=
\left(
\left(
\begin{array}{cc}
1 & 0 \\
0 & 1
\end{array}
\right)
+
\Psi_{-1}^{\infty}
\lambda^{-1}
+
\cdots
\right)
\exp T^{\infty}(\lambda), \\
\left(\Psi_{-1}^{\infty}\right)_{11}
&=
\sum_{k=1}^g
\frac{\partial}{\partial z_k}
\log \Theta [\p,\q](\mathbf{v}(t))
\frac{dU_k}{dx_{\infty^1}}(\infty^1)
+
\frac{t}{2}
\left(
\frac16
S(\infty^1)-\frac{W(\infty^1,\infty^2)}{dx_{\infty^1}dx_{\infty^2}}
\right),
\\
\left(\Psi_{-1}^{\infty}\right)_{21}
&=
\frac{i}{E(\infty^2,\infty^1)}
\frac{\Theta[\p,\q](2U(\infty^1)+\mathbf{v}(t))}{\Theta[\p,\q](\mathbf{v}(t))}
\exp
\left\{
t
\left(
c_{\infty^1}-c_{\infty^2}
\right)
\right\}, \\
\left(\Psi_{-1}^{\infty}\right)_{12}
&=
\frac{i}{E(\infty^1,\infty^2)}
\frac
{\Theta[\p,\q](2U(\infty^2)+\mathbf{v}(t))}
{\Theta[\p,\q](\mathbf{v}(t))}
\exp
\left\{
-t
\left(
c_{\infty^1}-c_{\infty^2}
\right)
\right\}, \\
\left(\Psi_{-1}^{\infty}\right)_{22}
&=
\sum_{k=1}^g
\frac{\partial}{\partial z_k}
\log \Theta[\p,\q](\mathbf{v}(t))
\frac{dU_k}{dx_{\infty^2}}(\infty^2)
-
\frac{t}{2}
\left(
\frac16 S(\infty^1)-\frac{W(\infty^1,\infty^2)}{dx_{\infty^1}dx_{\infty^2}}
\right).
\end{align*}
\end{lemma}

In the following theorem, 
we determine 
the monodromy matrices 
and the Stokes matrices of 
$\Psi(\lambda)$. 

\begin{theorem}
\label{thm:mono}
For $\nu = 1, 2, \ldots, 2g+2,$ 
the monodromy matrix $M_{\nu}$ of $\Psi(\lambda)$ 
corresponding to the contour 
$
l_{\nu} 
$ 
is given by
\begin{equation}
M_{\nu} 
=
\left(
\begin{array}{cc}
0                 & m_{\nu} \\
- m_{\nu}^{-1} & 0
\end{array}
\right),
\end{equation}
where
\begin{align*}
&  
m_{1} = -i, \,\, 
m_{2} = i (-1)^{g+1}\exp \left\{-2 \pi i \sum_{k=1}^g p_k \right\},  \\
&  
m_{2j+1} = i (-1)^g  \exp \left\{ 2 \pi i q_j-2\pi i\sum_{k=j}^g p_k \right\}, \\
&  
m_{2j+2} = i (-1)^{g+1} \exp \left\{ 2\pi i q_j -2\pi i \sum_{k=j+1}^g p_k \right\},
\end{align*}
for $j=1,2,\ldots,g$. 
The Stokes matrices are expressed by
$$
S_1^{\infty} 
=
S_2^{\infty}
=
\left(
\begin{array}{cc}
1 & 0  \\
0 & 1 \\
\end{array}
\right).
$$
\end{theorem}

\begin{proof}
By the involution $*$, we get
$$
\Psi(\lambda) M_{1} 
=
\Psi(\lambda) 
\left(
\begin{array}{cc}
  & i \\
i &     
\end{array}
\right),
\,
\textrm{or}
\,
=
\Psi(\lambda) 
\left(
\begin{array}{cc}
   & -i \\
-i &     
\end{array}
\right).
$$
We define 
\begin{align*}
M_{1}
=
\left(
\begin{array}{cc}
                  & m_{1}  \\
-m_{1}^{-1}  &
\end{array}
\right)
=
\left(
\begin{array}{cc}
    & -i  \\
-i  &
\end{array}
\right).
\end{align*}
Proposition \ref{prop:reg} implies that
\begin{align*}
T_{a_j}
\left[
\Psi(\lambda)
\right]
&=
\Psi(\lambda) M_{2j+2} M_{2j+1}  \\
&=
\Psi(\lambda)
\frac{T_{l_{2j+1}\circ l_{2j+2}}[\sqrt{\det \Phi(P)}]}{\sqrt{\det \Phi(P)}}
\exp \{2 \pi i (p_j+p_j^S) \}
\sigma_3,
\end{align*} 
and
\begin{align*}
T_{-b_j+b_{j-1}}
\left[
\Psi(\lambda)
\right]
&=
\Psi(\lambda) M_{2j+1} M_{2j}  \\
&=
\Psi(\lambda)
\frac{T_{l_{2j}\circ l_{2j+1}}[\sqrt{\det \Phi(P)}]}{\sqrt{\det \Phi(P)}}
\exp \{2 \pi i (q_j-q_{j-1}+q_j-q_{j-1})\}
\sigma_3,
\end{align*}
where 
$
\sigma_3=
\left(
\begin{array}{cc}
1 & 0 \\
0 &-1
\end{array}
\right).
$ 
In order to determine the monodromy matrices, 
we have 
\begin{alignat*}{3}
&U(\lambda_1)=0, & 
&U(\lambda_2) \,\, = \,\, \frac12 \sum_{k=1}^g \e_k, \\
&U(\lambda_{2j+1})=\frac12 \B \e_j +\frac12 \sum_{k=j}^g \e_k, & \quad
&U(\lambda_{2j+2}) \, = \, \frac12 \B \e_j +\frac12 \sum_{k=j+1}^g \e_k, \,\,
(j=1,2,\ldots,g).
\end{alignat*}
Then, 
we get 
\begin{equation*}
p_j^S=\frac12 
\left(
\delta_{2j+1}+\delta_{2j+2}+1
\right), \,
q_{j+1}-q_{j}
=
\frac12
\left(
\delta_{2j+2}+\delta_{2j+3}+1
\right),
\end{equation*}
where
\begin{equation*}
\begin{cases}
\delta_j=1 \quad \mathrm{if} \,\, \lambda_j \in S \\
\delta_j=0 \quad \mathrm{if} \,\, \lambda_j \notin S.
\end{cases}
\end{equation*} 
\par
The function 
$\sqrt{\det \Phi(P)}$ transforms 
with respect to the tracing along the cycles $l_j$ 
in the following way:
\begin{equation}
\label{eqn:root}
\begin{cases}
T_{l_{2j+1}\circ l_{2j+2}}[\sqrt{\det \Phi(P)}]
=
\exp \{\pi i (\delta_{2j+1}+\delta_{2j+2}+1) \}
\sqrt{\det \Phi(P)}, \\
T_{l_{2j}\circ l_{2j+1}}[\sqrt{\det \Phi(P)}]
=
\exp \{\pi i (\delta_{2j+2}+\delta_{2j+3}+1) \}
\sqrt{\det \Phi(P)}.
\end{cases}
\end{equation}
In order to prove (\ref{eqn:root}), 
we have only to 
note that 
if $\lambda_j$ is in $S$, 
from Proposition \ref{prop:reg}, 
$\det \Phi(P)$ 
has a zero of order one 
at $\lambda=\lambda_j$ 
and 
that 
if $\lambda_j$ is not in $S$, 
from Proposition \ref{prop:reg}, 
$\det \Phi(P)$ 
has a zero of order three 
at $\lambda=\lambda_j$. 
Therefore, 
we obtain 
\begin{equation*}
\begin{cases}
M_{2j+2}M_{2j+1}=
\exp \{ 2 \pi i \sigma_3 \} \\
M_{2j+1} M_{2j} =
\exp \{2 \pi i (q_j-q_{j-1}) \sigma_3 \}.
\end{cases}
\end{equation*} 
By 
considering 
\begin{equation*}
M_{2g+2} M_{2g+1}
\cdots 
M_{1}
=
I,
\end{equation*} 
we get
the monodromy matrices. 
\newline
\par
Because $\Psi(\lambda) \exp\left\{-T^{\infty}(\lambda)\right\}$ 
can be developed near $\lambda=\infty$ as a convergent series, 
the Stokes multipliers are 
$
\left(
\begin{array}{cc}
1 & 0 \\
0 & 1
\end{array}
\right).
$
\end{proof}

We can describe the monodromy data of $\Psi(\lambda)$ 
in the following way.

\begin{corollary}
$\Psi(\lambda)$ has the following monodromy data:
\newline
\newline
$\mathrm{(i)}$ \quad Stokes multipliers $S^{\infty}_1=S^{\infty}_2=
\left(
\begin{array}{cc}
1 & 0 \\
0 & 1
\end{array}
\right)
,$
\newline
$\mathrm{(ii)}$ \quad connection matrices 
$
C_j
=
\frac{1}{\sqrt{2 i m_{j}}}
\left(
\begin{array}{cc}
1 &  i m_{j} \\
-1 &  i m_{j} 
\end{array}
\right)
\,\,
(j=1,2, \ldots, 2g+2),
$
\newline
$\mathrm{(iii)}$ \quad exponents of formal monodromy 
$
T_{j,0}
=
\diag(-\frac14, \frac14)
\,\,
(j=1,2,\ldots,2g+2).
$
\newline
\newline
Especially, 
the developments of $\Psi(\lambda)$ 
near 
$
\lambda=\lambda_{j} \, 
(j=1,2,\ldots, 2g+2)
$ 
are expressed by
\begin{align*}
& 
\Psi(\lambda)
=
G_{j}
(
1
+
O(\lambda-\lambda_{j})
)
\exp T_{j} (\lambda) 
C_{j}, \\
& 
T_{j}(\lambda)
=
T_{j,0} \log (\lambda-\lambda_j ).
\end{align*}
\end{corollary}

\begin{proof}
(i) is clear. 
(ii) and (iii) can be obtained 
by diagonalizing the monodromy matrices 
$M_{j} \,\, (j=1,2,\ldots,2g+2).$
\end{proof}

This corollary 
means 
that 
the monodromy data of $\Psi(\lambda)$ 
are independent of 
the deformation parameters $t, \lambda_1, \lambda_2, \ldots, \lambda_{2g+2}.$

\section{Monodromy-Preserving Deformation}
In this section, 
we prove that 
$\Psi(\lambda)$ satisfies an ordinary differential 
equation which is expressed by
\begin{equation*}
\frac{d \Psi}{d \lambda}
=
\left(
\sum_{j=1}^{2g+2}
\frac{A_{j}}{\lambda-\lambda_{j}}
-
B_{-1}
\right)
\Psi(\lambda),
\end{equation*}
and 
concretely 
determine 
the coefficients 
$B_{-1}, A_{j} \,\, (1 \leq j \leq 2g+2).$
\par
We note 
that 
the monodromy matrices 
and 
the Stokes matrices 
are independent of the parameters $P_{\varphi}, P_{\phi}, p_k^S, q_k^S\,\,(1\leq k \leq g)$ 
in Theorem \ref{thm:mono}.
Thus, in this section, 
we set 
$
P_{\varphi}=\infty^2, 
\, 
P_{\phi}=\infty^1
$ 
and 
take $S_j$ 
so that $\lambda_j$ is not in $S_j$ 
because of the uniqueness of a solution of the Riemann-Hilbert problem. 
Therefore, 
we get
\begin{equation*}
\Psi(\lambda)
=
\frac{1}{\sqrt{\det \Phi^{\infty}(\lambda)}}
\Phi^{\infty}(\lambda),
\end{equation*}
where 
\begin{align*}
& 
\Psi(P)
=
\left(
\begin{array}{cc}
\varphi_j^{\infty}(P) & \varphi_j^{\infty} (P^{*}) \\
\psi_j^{\infty}(P)    & \psi_j^{\infty} (P^{*}) 
\end{array}
\right), \\
& 
\varphi_j^{\infty}(P)=
\frac
{
\Theta[\p,\q](U(P)+U(\infty^2)+\mathbf{v}(t))
\Theta[\p^{S_j},\q^{S_j}](U(P)-U(\infty^2))
}
{
\Theta[\p,\q](\mathbf{v}(t))
\Theta[\p^S,\q^S](-2U(\infty^2))
} \\
& 
\hspace{30mm}
\times
\exp \{-\frac{t}{2}(c_{\infty^1}-c_{\infty^2}) \}
\exp \Pi (P),
\\
& 
\psi_j^{\infty}(P)=
\frac
{
\Theta[\p,\q](U(P)+U(\infty^1)+\mathbf{v}(t))
\Theta[\p^{S_j},\q^{S_j}](U(P)-U(\infty^1))
}
{
\Theta[\p,\q](\mathbf{v}(t))
\Theta[\p^S,\q^S](-2U(\infty^1))
} \\
& 
\hspace{30mm}
\times
\exp \{\frac{t}{2}(c_{\infty^1}-c_{\infty^2}) \}
\exp \Pi (P).
\end{align*}

\begin{theorem}
\label{thm:mpd}
$\Psi(\lambda)$ satisfies 
the following ordinary differential equation:
\begin{equation}
\label{eqn:mpd}
\frac{d \Psi}{d \lambda}
= 
\left( 
\sum_{j=1}^{2g+2}
\frac{A_{j}}{\lambda-\lambda_{j}} 
-
B_{-1}
\right)
\Psi(\lambda),
\end{equation}
where
\begin{align*}
& 
B_{-1}
=
\diag (\frac{t}{2},-\frac{t}{2} ), \\
& 
A_j 
=
-\frac14
F_j^{\infty} \sigma_3 (F_j^{\infty})^{-1}, \quad 
\sigma_3=    
\left(
\begin{array}{cc}
1 & 0 \\
0 & -1
\end{array}
\right),
\\
& 
F_j^{\infty} 
=
\left(
\begin{array}{cc}
\varphi_j^{\infty}(\lambda_j) & \frac{d}{dx_j} \varphi_j^{\infty}(\lambda_j) \\   
\psi_j^{\infty}(\lambda_j) & \frac{d}{dx_j} \psi_j^{\infty}(\lambda_j)   
\end{array}
\right), 
\end{align*}
and $x_j=\sqrt{\lambda-\lambda_j}.$
\end{theorem}

\begin{proof}
From Lemma \ref{lem:ya}, it follows that
\begin{equation*}
\Psi_{\lambda} \Psi^{-1}
=
\left(
\begin{array}{cc}
-\frac{t}{2} &             \\
             & \frac{t}{2} 
\end{array}
\right)
+
O(\lambda^{-1})
:=
-B_{-1}+O(\lambda^{-1}).
\end{equation*}
\par
Near $\lambda=\lambda_j,$ 
we have
\begin{equation*}
\begin{cases}
\varphi_j^{\infty}(P)
=
\varphi_j(\lambda_j)
+
\sqrt{\lambda-\lambda_j}
\frac{d}{dx_j}\varphi_j(\lambda_j)
+\cdots \\
\psi_j^{\infty}(P)
=
\psi_j(\lambda_j)
+
\sqrt{\lambda-\lambda_j}
\frac{d}{dx_j}\psi_j(\lambda_j)
+\cdots, 
\end{cases}
\end{equation*}
which implies that
\begin{equation*}
\det \Phi^{\infty}(P)
=
-2
\sqrt{\lambda-\lambda_j}
\det F_j^{\infty}
+
O(\lambda-\lambda_j).
\end{equation*}
From the definition of $S_j$, 
it follows that $\det F^{\infty}_j \neq 0.$ 
\par
We set 
\begin{equation*}
\Psi(\lambda)
:=
G_j
(1+O(\lambda-\lambda_j))
(\lambda-\lambda_j)^
{
\left(
\begin{array}{cc}
-\frac14 &   \\
         & \frac14
\end{array}
\right)
}
C_j.
\end{equation*}
Then, 
we get
\begin{align*}
G_j=
\left.
\Psi(\lambda)
C_j^{-1}
(\lambda-\lambda_j)^
{
\left(
\begin{array}{cc}
\frac14 &   \\
         & -\frac14
\end{array}
\right)
}
\right|_{\lambda=\lambda_j}
&=
(-2\det F^{\infty}_j)^{-\frac12}
(1+O(\sqrt{\lambda-\lambda_j}))  \\
& 
\hspace{10mm}
\left.
\times
\Phi^{\infty}(P)
C_j^{-1}
(\lambda-\lambda_j)^
{
\left(
\begin{array}{cc}
1 &   \\
         & -\frac12
\end{array}
\right)
} 
\right|_{\lambda=\lambda_j} \\
&=
(-2 \det F^{\infty}_j)^{-\frac12}
\times
2i m_j
F^{\infty}_j.
\end{align*}
Thus, 
we obtain 
\begin{equation*}
\Psi_{\lambda}\Psi^{-1}
=
-\frac14
\frac{1}{\lambda-\lambda_j}
F^{\infty}_j \sigma_3 (F^{\infty}_j)^{-1}
+
\cdots
=
\frac{A_j}{\lambda-\lambda_j}
+\cdots.
\end{equation*}
\end{proof}

From Theorem \ref{thm:mpd}, 
we can obtain the following deformation equation:

\begin{corollary}
The deformation equation 
of the monodromy-preserving deformation (\ref{eqn:mpd}) 
is as follows. For $j,k=1,2,\ldots,2g+2,$
\begin{align*}
&
dA_j =   [\Theta_j,A_j], \\
&
d F^{\infty}_j = \Theta_j A_j, \\
&
\Theta_j =
\sum_{k \neq j}
A_k 
\frac{d\lambda_k-d\lambda_j}{\lambda_k-\lambda_j}
-
\left[
\Psi^{\infty}_{-1},dB_{-1}
\right]
-
d
\left(
\lambda_j B_{-1}
\right),
\end{align*}
where $``\,d\,"$ is the exterior differentiation 
with respect to the deformation parameters, 
$t, \lambda_1, \lambda_2, \ldots \lambda_{2g+2}$. 
\end{corollary}

\begin{proof}
See \cite{MM1}.
\end{proof}

\section{The $\tau$-Function for the Schlesinger System}
In this section, 
we calculate 
the $\tau$-function for the monodromy-preserving deformation 
(\ref{eqn:mpd}). 
This section 
consists of 
three subsections. 
Subsection 5.1 is devoted to 
the Hamiltonian $H_t$. 
Subsection 5.2, 5.3 
is devoted to the Hamiltonian 
$H_{j} \, (j=1,2,\ldots, 2g+2).$ 
In subsection 5.2, 
we quote Fay's identities and Rauch's variational formulas 
in order to compute the $\tau$-function. 
In subsection 5.3, 
we calculate $H_{j}$ and the $\tau$-function.

\subsection{The Hamiltonian at the Irregular Singular Point}
In this subsection, 
we prove Proposition \ref{prop:omegaa}, 
where we compute $\omega_{\infty}$ and the Hamiltonian $H_t$. 
\begin{proposition}
\label{prop:omegaa}
\begin{align*}
\omega_{\infty} 
&=
\left( 
\frac{1}{\Theta [\p,\q] ({\it \mathbf{v}}(t))}
\sum_{k=1}^g 
\frac{\partial}{\partial z_k} 
\{\Theta [\p,\q](\mathbf{v}(t)) \}
\frac{dU_k}{dx_{\infty^1}}(\infty^1)
+
\frac{t}{2} 
\left( 
\frac16 S(\infty^1)-
\frac{W(\infty^1,\infty^2)}{dx_{\infty^1}dx_{\infty^2}}
\right)
\right) dt \\
&= H_t dt.
\end{align*}

\end{proposition}

\begin{proof}
We define 
\begin{equation}
\Pi(P)
=
-\frac{t}{2} 
\{ 
-
\int^P_{\lambda_1}
W(P,\infty^1)
+
\int^P_{\lambda_1}
W(P,\infty^2)
\},  
\end{equation}
and let $\hat{\Pi}(P)$ denote the regular part of $\Pi(P)$ 
around $\lambda=\infty^1$ which is given by
\begin{equation}
\hat{\Pi}(P) 
= 
-\frac{t}{2}
\{
const \, + 
\left(-\frac16 S(\infty^1)+\frac{W(\infty^1,\infty^2)}{dx_{\infty^1}dx_{\infty^2}}\right)\lambda^{-1}
+\cdots 
\}.
\end{equation} 
Furthermore, 
we set 
\begin{align*}
\hat{\varphi}(P)
&=
\Theta[\p,\q](U(P)+U(P_{\varphi})+\mathbf{v}(t))
\Theta [\p^S,\q^S](U(P)-U(P_{\varphi})), \\
\hat{\psi}(P)
&=
\Theta[\p,\q](U(P)+U(P_{\psi})+\mathbf{v}(t))
\Theta[\p^S,\q^S](U(P)-U(P_{\psi})).
\end{align*}
Then, we get
\begin{align*}
\Psi(P) 
&=
\frac
{
\sqrt{\det \Phi(\infty)}
}
{
\sqrt{\det \Phi(P)}
} 
(G^{\infty})^{-1} 
\Phi(P) \\
&=
\frac
{
\sqrt{\det \Phi(\infty)}
}
{
\sqrt{\det \Phi(P)}
} 
(G^{\infty})^{-1}
\left(
\begin{array}{cc}
\hat{\varphi}(P) \exp (\hat{\Pi} (P)) 
& 
\hat{\varphi}(P^{*}) \exp (\hat{\Pi} (P^{*})) \\
\hat{\psi}(P) \exp (\hat{\Pi} (P)) 
& 
\hat{\psi}(P^{*}) \exp (\hat{\Pi} (P^{*}))
\end{array}
\right) \\
& 
\hspace{50mm}
\times
\diag 
\left(
\exp
\Big\{
- 
\frac{t}{2} 
\lambda
\Big\}, 
\exp
\Big\{
\frac{t}{2} 
\lambda
\Big\}
\right) \\
&:=
\hat{\Psi}^{\infty}(\lambda) 
\exp T^{\infty} (\lambda).
\end{align*}
\par
From the definition of $\omega_{\infty},$ 
it follows that
\begin{equation}
\omega_{\infty} 
= 
- 
\Res_{\lambda=\infty} 
\trace \hat{\Psi}^{\infty}(\lambda)^{-1} 
\frac{\partial}{\partial \lambda} 
\hat{\Psi}^{\infty}(\lambda) 
d T^{\infty} (\lambda).
\end{equation}
In order to compute $\omega_{\infty},$ we set 
\begin{equation}
A(\lambda) 
= 
(G^{\infty})^{-1}
\left(
\begin{array}{cc}
\hat{\varphi}(P) \exp \hat{\Pi} (P) 
& 
\hat{\varphi}(P^*) \exp \hat{\Pi}(P^*) \\
\hat{\psi} (P) \exp \hat{\Pi} (P)   
& \hat{\psi}(P^*) \exp \hat{\Pi}(P^*)
\end{array}
\right).
\end{equation}
Therefore, we get
\begin{equation}
- 
\trace 
\hat{\Psi}^{\infty}(\lambda)^{-1} 
\frac{\partial}{\partial \lambda} 
\hat{\Psi}^{\infty} (\lambda) 
d T^{\infty} (\lambda)  
=
\frac
{d t}
{2}  
\lambda 
\trace A^{-1}(\lambda) A^{\prime}(\lambda) 
\left(
\begin{array}{cc}
1 &  \\
  & -1 
\end{array}
\right),
\end{equation}
and
\begin{align}
\label{eqn:trace1}
&
\trace A^{-1}(\lambda) A^{\prime}(\lambda) 
\left(
\begin{array}{cc}
1 &  \\
  & -1 
\end{array}
\right)  \notag \\
&= 
\frac{1}{\det \Phi (P)}
\Big[
\det
\left(
\begin{array}{cc}
\{ 
\hat{\varphi} (P) \exp \hat{\Pi}(P) 
\}^{\prime} 
& \hat{\varphi}(P^*) \exp \hat{\Pi}(P^*) \notag \\
\{ 
\hat{\psi} (P) \exp \hat{\Pi}(P) 
\}^{\prime}    
& 
\hat{\psi}(P^*) \exp \hat{\Pi}(P^*) 
\end{array}
\right) \\      
&  
\hspace{40mm}           
-
\det
\left(
\begin{array}{cc}
\hat{\varphi} (P) \exp \hat{\Pi}(P) 
& 
\{ 
\hat{\varphi}(P^*) \exp \hat{\Pi}(P^*) 
\}^{\prime} \\
\hat{\psi} (P) \exp \hat{\Pi} (P) 
& 
\{ \hat{\psi}(P^*) \exp \hat{\Pi} (P^*) \}^{\prime}
\end{array}
\right)
\Big],
\end{align}
where $\prime$ means the differentiation 
with respect to the variable $\lambda$.
\par
We have normalized the matrix function $\Psi(\lambda)$ around $\lambda=\infty$ 
in Lemma \ref{lem:ya} 
and have proved that 
the monodromy data of $\Psi(\lambda)$ 
are independent of $P_{\varphi}, P_{\psi}$ 
in Theorem \ref{thm:mono} and its corollary. 
Therefore, 
we can choose the parameters, 
$P_{\varphi}, P_{\psi}$ 
at our disposal to simplify the calculation. 
\par
Firstly, we multiply 
both the numerators and the denominators of (\ref{eqn:trace1}) 
by $\displaystyle\frac{1}{\lambda_{\psi}-\lambda_{\varphi}}$. 
Then, 
we take the limit $P_{\psi} \rightarrow P_{\varphi}$ 
and get
\begin{equation}
\hat{\psi}(P) 
= 
\frac
{
\partial \hat{\varphi}(P)
}
{\partial \lambda_{\varphi}}.
\end{equation}
Next, we multiply 
both 
the numerators 
and 
the denominators of (\ref{eqn:trace1}) 
by $\displaystyle\frac{1}{\lambda_{\varphi}-\lambda}$. 
Then, 
we take the limit 
$
P_{\varphi} \rightarrow P
$ 
and obtain
\begin{align}
\trace A^{-1}(\lambda) A^{\prime}(\lambda)
\left(
\begin{array}{cc}
1 &  \\
  & -1 
\end{array}
\right)
&=
2 
\frac{1}{\Theta[\p,\q](\mathbf{v}(t))} 
\lim_{P_{\varphi} \rightarrow P} 
\frac{\partial}{\partial x_{\varphi}}
\Theta[\p,\q](-U(P) + U(P_{\varphi}) + \mathbf{v}(t))  \notag \\
& 
\hspace{50mm}
+ 
2 
\frac{\partial}{\partial \lambda}
\{
\hat{\Pi} (P)
\}.
\end{align}
From the definition of $\omega_{\infty},$ 
it follows that
\begin{equation*}
\omega_{\infty} = 
\frac{dt}{\Theta [\p,\q] (\mathbf{v}(t))}
\sum_{k=1}^g 
\frac{\partial}{\partial z_k} 
\{\Theta [\p,\q](\mathbf{v}(t)) \}
\frac{dU_k}{dx_{\infty^1}}(\infty^1)
+
\frac{t}{2}dt 
\left( 
\frac16 S(\infty^1)-\frac{W(\infty^1,\infty^2)}{dx_{\infty^1}dx_{\infty^2}}
\right).
\end{equation*}
\end{proof}

\subsection{Fay's Identities and Rauch's Variational Formulas}
In this subsection, 
we quote 
Fay's identities and Rauch's variational formulas 
in order to determine the Hamiltonians $H_j\,(j=1,2,\ldots,2g+2).$ 
Lemma \ref{lem:Fay2} and \ref{lem:Rauch} 
are devoted to 
Fay's Identities and Rauch's Variational Formulas, 
respectively.

\begin{lemma}
\label{lem:Fay2}
{\it
(1) 
\quad
For $P,Q \in \mathcal{L},$
\begin{equation*}
\frac{
\Theta [\mathbf{p}^{T}, \mathbf{q}^{T} ]^2(U(P)-U(Q))}
{
\Theta [\mathbf{p}^{T}, \mathbf{q}^{T} ]^2(0)
E^2(P,Q)
}
=W(P,Q)+
\sum_{k,l=1}^g 
\frac{\partial^2}{\partial z_k \partial z_l}
\log \Theta [\mathbf{p}^{T}, \mathbf{q}^{T} ](0)
dU_k(P) dU_l(Q).
\end{equation*}
(2) \quad 
For $P,Q \in \mathcal{L},$
\begin{align*}
&
\frac{\Theta [\mathbf{p}^{T}, \mathbf{q}^{T} ](2(U(P)-U(Q)))}
{
\Theta [\mathbf{p}^{T}, \mathbf{q}^{T} ](0)
E^4(P,Q)
}
-
\frac{\Theta [\mathbf{p}^{T}, \mathbf{q}^{T} ]^4(U(P)-U(Q))}
{
\Theta [\mathbf{p}^{T}, \mathbf{q}^{T} ]^4(0)
E^4(P,Q)
} \\
& 
\hspace{30mm}
=
\frac12 
\sum_{k,l,m,n=1}^g 
\frac{\partial^4}{\partial z_k \partial z_l \partial z_m \partial z_n}
\log \Theta [\mathbf{p}^{T}, \mathbf{q}^{T} ](0) 
dU_k(P) dU_l(P)dU_m(Q) dU_n(Q).
\end{align*}
(3) \quad 
For $P,Q\in\mathcal{L},$ 
\begin{align*}
\frac
{\Theta[\mathbf{p}^T, \mathbf{q}^T](0)\Theta[\mathbf{p},\mathbf{q}]\left(2\left(U(P)-U(Q)\right)\right)}
{\Theta[\mathbf{p}^T,\mathbf{q}^T]\left(U(P)-(Q)\right)^2E(P,Q)^2}
&=
W(P,Q)               \\
&\quad\quad
+
\sum_{k,l=1}^g
\frac{\partial^2}{\partial z_k \partial z_l}
\log \Theta[\mathbf{p}^T,\mathbf{q}^T]\left(U(P)-U(Q)\right)dU_k(P)dU_l(Q) \\
&=
dx_{P}dx_{Q}
\log
\frac{E(P,Q)}{\Theta[\mathbf{p}^T,\mathbf{q}^T]\left(U(P)-U(Q)\right)},
\end{align*}
where $x_P,x_Q$ are local coordinates of $P,Q,$ respectively.
Especially, since $\mathcal{L}$ is hyperelliptic, 
it follows that 
\begin{equation*}
\frac
{\Theta[\mathbf{p}^T, \mathbf{q}^T](0)\Theta[\mathbf{p},\mathbf{q}]\left(2\left(U(P)-U(Q)\right)\right)}
{\Theta[\mathbf{p}^T,\mathbf{q}^T]\left(U(P)-(Q)\right)^2E(P,Q)^2}
= 
dx_{P} dx_{Q}\log \frac{1}{m_T(P,Q)}.
\end{equation*}
}
\end{lemma}

\begin{proof}
For (1), (2) and (3),  
see pp 26, pp 28 and pp 29 in Fay's book \cite{Fay}, 
respectively.
\end{proof}

\par
Rauch \cite{Rauch} described 
the dependence of $dU_k\,(k=1,2,\ldots,g)$ and $\B_{kl}\,(k,l=1,2,\ldots,g)$ 
on the moduli of the Riemann surfaces. 
The moduli space of hyperelliptic curves can be parameterized 
by the positions of the branch points $\lambda_j \, (j=1,2,\ldots,2g+2)$. 
Korotkin \cite{Korotkin} proved the variational formulas of the following useful form.

\begin{lemma}
\label{lem:Rauch}
{\it
(1) \quad 
For $P\in \mathcal{L},$ 
\begin{equation}
\label{eqn:Rauch1}
\frac{\partial}{\partial \lambda_j} 
\left\{
\frac{dU_k}{dx_{P}}(P)
\right\}
=
\frac12 \frac{W(P,\lambda_j)}{dx_{P} d x_j}
\frac{dU_k}{dx_j}(\lambda_j),
\end{equation}
where $x_{P}$ is a local coordinate of $P$ and 
$x_j=\sqrt{\lambda-\lambda_j},$ which 
is a local coordinate of the branch point $\lambda_j$ 
for any $j=1,2,\ldots, 2g+2.$
\newline
(2) \quad 
For the branch points $\lambda_j \,\,(j=1,2,\ldots, 2g+2),$ 
\begin{equation}
\label{eqn:Rauch2}
\frac{\partial \B_{kl}}{\partial \lambda_j} 
=
\pi i 
\frac{dU_l}{dx_j}(\lambda_j)
\frac{dU_k}{dx_j}(\lambda_j) \quad (k,l=1,2,\ldots, g),
\end{equation}
where $x_j=\sqrt{\lambda-\lambda_j},$ which 
is a local coordinate of $\lambda_j.$
}
\end{lemma}

\subsection{The $\tau$-Function}
In this subsection, 
we compute $\omega_{\lambda_j} (j=1,2,\ldots,2g+2)$ and 
calculate the $\tau$-function.

\begin{proposition}
\label{lem:omeganu}
For $j=1,2,\ldots,2g+2,$ 
\begin{align*}
H_j
&=
\frac{\partial}{\partial \lambda_j}
\log \Theta[\p,\q](\mathbf{v}(t) | \B)
-\frac12
\frac{\partial}{\partial \lambda_j}
\log \det \mathcal{A} 
-\frac18
\frac{\partial}{\partial \lambda_j}
\log 
\prod_{k<l} 
(\lambda_k-\lambda_l) \\
&
\hspace{60mm}
+ 
\frac{\partial}{\partial \lambda_j}
\left\{
\frac{t^2}{4}\left(\frac16 S(\infty^1)-\frac{W(\infty^1,\infty^2)}{dx_{\infty^1}dx_{\infty^2}}\right)
\right\}.
\end{align*}
\end{proposition}

\begin{proof}
By direct calculation, we obtain 
\begin{equation*}
\frac12 
\trace \left( \Psi^{\prime}(\lambda)\Psi^{-1}(\lambda) \right)^2
=
-
\frac{\det \left( \Phi_x \right)}{\det \Phi}
+
\frac14
\left(
\frac
{
(\det \Phi)_{\lambda}
}
{\det \Phi}
\right)^2.
\end{equation*}
We calculate $\omega_{\lambda_j}$ in the same way as $\omega_{\infty}$ 
in Proposition \ref{prop:omegaa}. 
\par
We multiply 
both the numerators and the denominators of
$
\displaystyle
\frac{\det \left( \Phi_{\lambda} \right)}{\det \Phi},
\frac
{
(\det \Phi)_{\lambda}
}
{\det \Phi}
$
by $\displaystyle\frac{1}{\lambda_{\varphi}-\lambda_{\psi}}$. 
Then, 
we take the limit $P_{\psi} \rightarrow P_{\varphi}$ 
and get
\begin{equation}
\psi(P) 
= 
\frac{\partial \varphi(P)}{\partial \lambda_{\varphi}}.
\end{equation}
Furthermore, 
we 
multiply 
both 
the numerators 
and 
the denominators of 
$
\displaystyle
\frac{\det \left( \Phi_{\lambda} \right)}{\det \Phi},
\frac
{
(\det \Phi)_{\lambda}
}
{\det \Phi}
$
by $\displaystyle\frac{1}{\lambda_{\varphi}-\lambda}$ 
and 
take the limit $P_{\varphi} \rightarrow P$. 
Then, we obtain 
\begin{align*}
\frac
{
(\det \Phi)_{\lambda}
}
{\det \Phi}
&=
2 
\frac{\partial}{\partial \lambda}
\log \Theta [\p^S,\q^S](-2U(P)),
\\
\frac{\det \left( \Phi_{\lambda} \right)}{\det \Phi} 
&= 
\frac{1}{\Theta[\p,\q](\mathbf{v}(t))}
\lim_{P_{\varphi} \rightarrow P}
\frac{\partial^2}{\partial \lambda \partial \lambda_{\varphi}}
\Theta[\p,\q](-U(P)+U(P_{\varphi})+\mathbf{v}(t))  \\
&+ 
\frac{2}{\Theta[\p,\q](\mathbf{v}(t))} 
\lim_{P_{\varphi} \rightarrow P}
\frac{\partial }{\partial \lambda}
\Theta[\p,\q] (-U(P)+U(P_{\varphi})+\mathbf{v}(t))
\frac{\partial }{\partial \lambda} \Pi (P) \\
&+
\frac{1}{\Theta[\p^S,\q^S](-2U(P))}
\lim_{P_{\varphi} \rightarrow P}
\frac{\partial^2}{\partial \lambda \partial \lambda_{\varphi}}
\Theta [\p^S,\q^S] (-U(P)-U(P_{\varphi})) \\
&-
\left(
\frac{\partial}{\partial \lambda} \Pi (P)
\right)^2,
\end{align*}
which implies that
\begin{align}
& 
\frac12 
\trace \left( \Psi^{\prime}(\lambda)\Psi^{-1}(\lambda) \right)^2
= 
-
\left.
\frac{\partial^2}{\partial \lambda \partial \lambda_{\varphi}}
\log \Theta [\p^S,\q^S] (-U(P) -U(P_{\varphi})) 
\right|_{P_{\varphi}=P}  \\
&
\hspace{25mm}
-  
\frac{1}{\Theta[\p,\q](\mathbf{v}(t))}
\left.
\frac{\partial^2}{\partial \lambda \partial \lambda_{\varphi}}
\Theta[\p,\q] (-U(P) + U(P_{\varphi}) +\mathbf{v}(t))
\right|_{P_{\varphi}=P}  \\
&
\hspace{25mm}
-  
\left.
\frac{2}{\Theta[\p,\q](\mathbf{v}(t))}
\frac{\partial}{\partial \lambda} \Theta[\p,\q](-U(P)+U(P_{\varphi})+\mathbf{v}(t))
\frac{\partial}{\partial \lambda}
\Pi (P)  
\right|_{P_{\varphi=P}}  \\
& 
\hspace{25mm}
+ 
\left(
\frac{\partial}{\partial \lambda}
\Pi (P)
\right)^2.
\end{align}
\par
Firstly, we calculate the residue of (5.12) at $\lambda=\lambda_j$. 
For the purpose, 
we set 
$$
P=\lambda_j, \,x_{P}=x_j:=\sqrt{\lambda-\lambda_j}
$$ 
in (\ref{eqn:proj2}) of Lemma \ref{lem:Fay1}.
Then, we have

\begin{align*}
&
\Res_{\lambda=\lambda_1} 
\left\{
\left.
-
\frac{\partial^2}{\partial \lambda \partial \lambda_{\varphi}}
\log \Theta [\p^S,\q^S] (-U(P) -U(P_{\varphi})) 
\right|_{P_{\varphi}=P}
\right\}
\\
&
\hspace{20mm}
=
\frac18 \sum_{k\neq j} \frac{n_j n_k}{\lambda_j-\lambda_k} 
-\frac{1}{4 \Theta [\p^T,\q^T]}(0)
\sum_{l,k=1}^g \frac{\partial^2 \Theta[\p^T,\q^T]}{\partial z_l \partial z_k}(0)
\frac{dU_l}{dx_j}(\lambda_j) \frac{dU_k}{dx_j}(\lambda_j),
\end{align*}
where 
$n_k=1$ for $\lambda_k \in T$ and $n_k =- 1$ for $\lambda_k \notin T.$ 
For the calculation, 
we use (\ref{eqn:Rauch2}) in Lemma \ref{lem:Rauch} and the heat equation 
\begin{equation}
\label{eqn:heat}
\frac{\partial^2 \Theta[\p,\q](\z|\B)}{\partial z_l \partial z_k}
=
4 \pi i 
\frac{\partial \Theta[\p,\q](\z|\B)}{\partial \B_{lk}}.
\end{equation}
Then, 
we get 
\begin{align*}
&
\Res_{\lambda=\lambda_j} 
\left\{
-
\left.
\frac{\partial^2}{\partial \lambda \partial \lambda_{\varphi}}
\log \Theta [\p^S,\q^S] (-U(P) -U(P_{\varphi})) 
\right|_{P_{\varphi}=P}
\right\} \\
&
\hspace{65mm}
=
\frac18 \sum_{k \neq j}\frac{n_j n_k}{\lambda_j-\lambda_k}
-
\frac{\partial}{\partial \lambda_j} 
\log 
\Theta [\p^T,\q^T](0 |\B).
\end{align*}
By using the Thomae's formula 
\begin{equation*}
\Theta^4[\p^T,\q^T](0)
=
\pm
\frac{(\det \mathcal{A})^2}{(2 \pi i)^{2g}}
\prod_{l<k,\,l,k=1}^{g+1} (\lambda_{i_l}-\lambda_{i_k})
\prod_{l<k,\,l,k=1}^{g+1} (\lambda_{j_l}-\lambda_{j_k}),
\end{equation*}
we get
\begin{align*}
&
\Res_{\lambda=\lambda_j} 
\left\{
-
\left.
\frac{\partial^2}{\partial \lambda \partial \lambda_{\varphi}}
\log \Theta [\p^S,\q^S] (-U(P) -U(P_{\varphi}))
\right|_{P_{\varphi}=P} 
\right\} \\
&
\hspace{60mm}
=
-\frac12
\frac{\partial}{\partial \lambda_j}
\log \det \mathcal{A} 
-\frac18
\frac{\partial}{\partial \lambda_j}
\log 
\prod_{k<l} 
(\lambda_k-\lambda_l). 
\end{align*} 
\par
Secondly, 
we calculate the residue of the sum of (5.13) and (5.14) at $\lambda=\lambda_j,$ 
which is
\begin{align}
\label{eqn:secandthir}
&
\frac{1}{4\Theta[\p,\q]( \mathbf{v}(t) )}
\sum_{k,l=1}^g
\frac{\partial^2 \Theta[\p,\q]}{\partial z_k \partial z_l}
(\mathbf{v}(t))
\frac{dU_k}{dx_j}(\lambda_j) \frac{dU_l}{dx_j}(\lambda_j)  \notag \\
&
\hspace{40mm}
+
\frac{1}{2\Theta[\p,\q](\mathbf{v}(t))}
\sum_{k=1}^g 
\frac{\partial \Theta[\p,\q]}{\partial z_k}
(\mathbf{v}(t))
\times t 
\frac{W(\lambda_1,\infty^1)}{dx_jdx_{\infty^1}}. 
\end{align}
From Lemma \ref{lem:Rauch}, it follows that 
(\ref{eqn:secandthir})  
is 
\begin{equation*}
\frac{\partial}{\partial \lambda_j} 
\log 
\Theta[\p,\q](\mathbf{v}(t)|\B).
\end{equation*}
\par 
Lastly, 
we deal with 
the residue of (5.15) at $\lambda=\lambda_j,$ 
which is 
$\displaystyle\frac{t^2}{4}\left(\frac{W(\lambda_j,\infty^1)}{dx_jdx_{\infty^1}}\right)^2.$ 
We prove that if $\lambda_j \in T,$
\begin{equation*}
\frac{\partial}{\partial \lambda_j}
\left\{
\frac{t^2}{4}
\left(
\frac16 S(\infty^1)-\frac{W(\infty^1,\infty^2)}{dx_{\infty^1}dx_{\infty^2}}
\right)
\right\}
=
\frac{t^2}{4}
\left(
\frac{
W(\lambda_j,\infty^1)
}{dx_jdx_{\infty^1}}
\right)^2.
\end{equation*}
If $\lambda_j \not\in T,$ 
this formula can be proved in the same way.  
\par
From Lemma \ref{lem:Fay1} and Lemma \ref{lem:Fay2} (1), 
it follows that 
\begin{align*}
\frac16 S(\infty^1)-\frac{W(\infty^1,\infty^2)}{dx_{\infty^1} dx_{\infty^2}}
&=
\frac{1}{16}
\left(
\sum_{k=1}^{g+1} 
\lambda_{i_k}-\sum_{k=1}^{g+1}\lambda_{j_k}
\right)^2
-
\frac
{
\Theta[\p^T,\q^T]^2(U(\infty^1)-U(\infty^2))
}
{
\Theta[\p^T,\q^T]^2(0)
\{E(\infty^1,\infty^2)\sqrt{dx_{\infty^1}}\sqrt{dx_{\infty^2}}\}^2
} \\
&\hspace{30mm}
-2
\sum_{k,l=1}^g 
\frac{\partial^2}{\partial z_k \partial z_l} 
\log \Theta [\p^T,\q^T](0) 
\frac{dU_k}{dx_{\infty^1}}(\infty^1)
\frac{dU_l}{dx_{\infty^1}}(\infty^1)  \\
&=
\frac18\left(\sum_{k=1}^{g+1} \lambda_{i_k}-\sum_{k=1}^{g+1}\lambda_{j_k}\right)^2  \\
&\hspace{30mm}
-2
\sum_{i,j=1}^g 
\frac{\partial^2}{\partial z_i \partial z_j} 
\log \Theta [\p^T,\q^T](0) 
\frac{dU_i}{dx_{\infty^1}}(\infty^1)
\frac{dU_j}{dx_{\infty^1}}(\infty^1), 
\end{align*}
where we have used 
\begin{equation*}
\frac
{
\Theta[\p^T,\q^T](U(\infty^1)-U(\infty^2))
}
{
\Theta[\p^T,\q^T](0)
E(\infty^1,\infty^2)
} 
=
m_T(\infty^1,\infty^2)
=
\frac{\sqrt{-1}}{4}
\left(\sum_{k=1}^{g+1} \lambda_{i_k}-\sum_{k=1}^{g+1}\lambda_{j_k}\right)\sqrt{dx_{\infty^1}}\sqrt{dx_{\infty^2}}.
\end{equation*}
Then, 
from 
Lemma \ref{lem:Rauch}, 
it follows that
\begin{align*}
&\frac{\partial}{\partial \lambda_j}
\left\{
\frac16 S(\infty^1)-\frac{W(\infty^1,\infty^2)}{dx_{\infty^1}dx_{\infty^2}}
\right\}  \\
&=
\frac14\left(\sum_{k=1}^{g+1} \lambda_{i_k}-\sum_{k=1}^{g+1}\lambda_{j_k}\right) \\
&
-
\frac12
\sum_{k,l,m,n=1}^g 
\frac{\partial^4}{\partial z_k \partial z_l \partial z_m \partial z_n}
\log \Theta[\mathbf{p}^T,\mathbf{q}^T](0) 
\frac{dU_k}{dx_{\infty^1}}(\infty^1)\frac{dU_l}{dx_{\infty^1}}(\infty^1)\frac{dU_m}{dx_{\infty^1}}(\infty^1)\frac{dU_n}{dx_{\infty^1}}(\infty^1)  \\
&
-
\left(
\sum_{k,l=1}^g \frac{\partial^2}{\partial z_k\partial z_l}
\log \Theta[\mathbf{p}^T,\mathbf{q}^T](0)\frac{dU_k}{dx_j}(\lambda_j)\frac{dU_l}{dx_{\infty^1}}(\infty^1)
\right)^2  \\
&
-2
\left(
\sum_{k,l=1}^g \frac{\partial^2}{\partial z_k\partial z_l}
\log \Theta[\mathbf{p}^T,\mathbf{q}^T](0)\frac{dU_k}{dx_j}(\lambda_j)\frac{dU_l}{dx_{\infty^1}}(\infty^1)
\right)
\frac{W(\lambda_j,\infty^1)}{dx_jdx_{\infty^1}}.
\end{align*}
Thus, from Lemma \ref{lem:Fay2} (1) and (2), 
we have 
\begin{align*}
\frac{\partial}{\partial \lambda_1}
\left\{
\frac16 S(\infty^1)-\frac{W(\infty^1,\infty^2)}{dx_{\infty^1}dx_{\infty^2}}
\right\}  
&=
\frac14\left(\sum_{k=1}^{g+1} \lambda_{i_k}-\sum_{k=1}^{g+1}\lambda_{j_k}\right) 
-
\frac{\Theta[\mathbf{p}^T,\mathbf{q}^T]\left(2\left(U(\infty^1)-U(\lambda_j)\right)\right)}
{\Theta[\mathbf{p}^T,\mathbf{q}^T](0)\{E(\lambda_j,\infty^1)\sqrt{dx_j}\sqrt{dx_{\infty^1}}\}^4}   \\
&+\left(\frac{W(\lambda_j,\infty^1)}{dx_jdx_{\infty^1}}\right)^2.
\end{align*}
Therefore, from Lemma \ref{lem:Fay2} (3), 
we obtain 
\begin{equation*}
\frac{\partial}{\partial \lambda_j}
\left\{
\frac16 S(\infty^1)-\frac{W(\infty^1,\infty^2)}{dx_{\infty^1}dx_{\infty^2}}
\right\}
=
\left(
\frac{W(\lambda_j,\infty^1)}{dx_jdx_{\infty^1}}\right)^2,
\end{equation*}
which completes the proof of the proposition.

\end{proof}

Since we have calculated the Hamiltonians 
$H_t, H_1, H_2, \ldots, H_{2g+2}$, 
we finally obtain Theorem 0.1.

\quad {\bf Acknowledgments.}
\, The author thanks 
Professor Yousuke Ohyama 
for the careful guidance 
and the referees for useful comments.

\end{document}